\documentclass[12pt,a4paper]{article}
\usepackage[T2A]{fontenc}
\usepackage[english]{babel}
\usepackage[utf8]{inputenc}

\usepackage{indentfirst}

\usepackage[hidelinks]{hyperref}

\usepackage{orcidlink}

\usepackage{amssymb}
\usepackage{changepage}
\usepackage{fdsymbol}
\usepackage{bussproofs}
\usepackage{xcolor}

\begin{document}

\renewcommand{\lq}{\guillemotleft}
\renewcommand{\rq}{\guillemotright}
\newcommand{\q}{\textrm{'}}
\renewcommand{\l}{\lceil}
\renewcommand{\r}{\rceil}
\newcommand{\la}{\langle}
\newcommand{\ra}{\rangle}
\newcommand{\sh}{{\sharp}}
\newcommand{\na}{{\natural}}
\renewcommand{\geq}{\geqslant}
\renewcommand{\leq}{\leqslant}
\newcommand{\name}{\mathsf{name}}
\newcommand{\arity}{\underline{\mathsf{arity}}}
\newcommand{\val}{\mathsf{val}}
\renewcommand{\S}{\mathsf{S}}
\renewcommand{\phi}{\varphi}

\title{\bfseries \itshape A Self-Reflecting Formal Language}

\author{Mikhail Patrakeev\footnote{E-mail: p17533@gmail.com; \href{https://orcid.org/0000-0001-7654-5208}{orcid.org/0000-0001-7654-5208}}\hspace{2mm}\orcidlink{0000-0001-7654-5208}
}

\providecommand{\keywords}[1]
{
  \small	
  \textbf{\textit{Keywords---}} #1
}

\maketitle

\begin{abstract}\rm
\begin{adjustwidth}{2cm}{2cm} 

We construct a formal theory, which we call reflectica, whose language possesses the following properties of natural language: it is a self-reflecting language and an intensional language. By a \textit{self-reflecting} language we understand an interpreted language that is a meta-language in relation to itself. 
By an \textit{intensional} language we understand a language that has expressive means sufficient to represent intensional features of a natural language, such as statements containing propositional attitude reports, various kinds of quotation, and other types of expressions with an intensional context. At the same time, we present a new method for constructing an intensional logic that allows us to make reflectica an intensional system much simpler than other well-known intensional logics.

\end{adjustwidth}
\end{abstract}


\section{Introduction}
\noindent

Some natural language expressions refer to (= denote, designate, connote) objects. But at the same time the words and the expressions of natural language are also objects that could be referred to by natural language expressions. This allows us to use natural language to describe and study that very natural language. For example, English allows linguists and philosophers to study English, in particular its syntax, semantics, and the relation that maps a natural language expression to the object it refers to. In this paper we construct a formal theory, which we call reflectica, whose language possesses the above property of natural language; we explicate this property as follows:

An \textit{interpretation} of a formal language ${L}$ is a partial\footnote{This means that the interpretation may be defined not on all expressions of $L$, but only on some of them.} function that maps the expressions of $L$ to their \textit{values} (= denotations, extensions, designations, connotations); we call the values of expressions of ${L}$ \textit{objects}, so that the set of objects equals by definition the range (= the image) of the interpretation of $L$. An \textit{interpreted language} is a language endowed with an interpretation. A \textit{self-interpreted language} is an interpreted language ${L}$ such that all symbols\footnote{Formal language symbols correspond to natural language words.} and expressions of $L$ are objects. We can use an interpreted language to express statements about object. Thus, if $L$ is a self-interpreted language, then we can use the language ${L}$ to express assertions about the symbols of ${L}$, the expressions of ${L}$, and the values of expressions of ${L}$ (because the latter are also objects). Therefore, if a self-interpreted language has sufficient expressive power, then we can use it to describe its own syntax, semantics, and interpretation; in which case we call such a language a \textit{self-reflecting language}.

Those familiar with axiomatic theories of truth~\cite{halbach2014axiomatic} may wonder: is not the language of a theory of truth a self-reflecting language? Most expositions of such theories do not build any interpretation; but to endow an axiomatic theory of truth with an interpretation that makes its language self-interpreted is quite a straightforward task. Also theories of truth are constructed in such a way that they allow us to formulate in their language the definition of truth that lies in the heart of their semantics; they also allow to describe their syntax. Thus, after adding appropriate interpretation, the language of such a theory will be almost self-reflecting: it will be a self-interpreted language that has expressive means to describe its own syntax and semantics. Is it possible to construct such an interpretation in a way that the interpreted language also allows us to describe the constructed interpretation? The methods of this paper can not give an answer to this question, since axiomatic theories of truth, being formulated in classical predicate logic, are extensional\footnote{An \label{footnote context} extensional theory is a theory without  intensional contexts. An intensional (also known as non-extensional or opaque) context is a part of a complex expression in which substitution of co-referring (i.e., having equal values) expressions may change the value of the whole complex expression. For example, the value of the expression \textit{Lois believes that Superman is Clark Kent} equals the \textit{false}, but it becomes the \textit{true} when we replace the subexpression \textit{Clark Kent} with the co-referring expression \textit{Superman}. Another example of an intensional context is a quoted expression like the subexpression \textit{one plus two} in the expression \textit{The number of words in the expression `one plus two'}.} theories, while reflectica is an intensional theory~\cite{sep-logic-intensional}. Therefore, a second question arises: Is it possible to construct an extensional formal theory whose language is a self-reflecting language? Note that if our goal is to construct a formal theory whose language is as close to natural language as possible, then the choice in favour of an intensional theory is an advantage rather than a disadvantage, since natural language is intensional.

The formal theory that we construct here, reflectica, has a self-reflecting language. Reflectica is an interpreted theory, it is endowed with an interpretation which is a partial function that maps reflectica expressions to their values. We call these values, including the truth-values, reflectica objects. The reflectica symbols and expressions are also objects. Most reflectica expressions contain the main symbol and several subexpressions. We consider the main symbol of an expression as an operation, the subexpressions as its arguments, and the value of the whole expression as the result of such an operation. In reflectica, every subexpression bears two objects: its value is an object and the subexpression itself is also an object. Therefore, the main symbol of a reflectica expression can act as an operation on objects that takes as input objects not only the values of the subexpressions, but also the subexpressions themselves. From the point of view of the use–mention distinction, we can express the same idea by saying that in reflectica a subexpression can be used (when its value is treated as an input object), can be mentioned (when the subexpression itself is treated as an input object), and can be both used and mentioned at the same time. 

Such a simple idea allows us to make reflectica an intensional system without introducing any additional entities. As a result, the expressive power of reflectica language is sufficient to represent intensional features of natural language, such as statements containing propositional attitude reports~\cite{sep-prop-attitude-reports} (like \textit{Lois believes that Superman is Clark Kent}), various kinds of quotation~\cite{sep-quotation} (such as \textit{Quine said that quotation `has a certain anomalous feature'}), and other types of natural language expressions with an intensional context, see Section~\ref{sect.context}. This new way of constructing an intensional logic makes reflectica much simpler than other well-known intensional systems~\cite{sep-logic-intensional}.

Reflectica allows us to write sentences that are propositions about reflectica, in particular, about its rules. For example, if we want to extend reflectica by adding a new rule, we can simply add a new axiom declaring that rule, see Section~\ref{sect.rules}. We can also write formulas with direct quantification over reflectica expressions, avoiding arith\-me\-ti\-za\-tion and G\"odel numbering, see Section~\ref{sect.quantifiers}. For example, the formula 
\[
\forall{\varphi}\,\big(\hspace{0.5pt}\mathsf{\underline{Form}}(\na\hspace{0.5pt}\varphi)\to[\varphi\:{\leftrightarrow}\:(\varphi=\top)]\big)
\]
is a reflectica sentence which says that every reflectica formula is equivalent to truthness of its value. In this sentence, the symbol $`{\na}\q$ is a quotation mark that does not bind variables, see Section~\ref{sect.quotes}; the symbols $`\mathsf{\underline{Form}}\q$ and $`{\to}\q$ can be found in Sections~\ref{sect.rules} and~\ref{sect.conditional}. Reflectica even allows us to correctly quantify over a variable that is inside a quotation (see the example at the end of Section~\ref{sect.context}), which some philosophers consider impossible, see BQ2 in~Section 2 of \cite{sep-quotation}.

\section{Brief informal exposition}
\noindent
Two important points:

\begin{itemize}
\item[1.] Every expression and every subexpression refers to (= denotes, designates) its value (= denotation, designation) whenever it has one.
\item[2.] Every object, including expressions and their values, is a finite string (= a sequence of symbols) in the alphabet of reflectica.
\end{itemize}

We denote (in the meta-theory) the value of expression ${e}$ by $\val({e})$ and use single quotes ` and ' to denote strings in the alphabet of reflectica. The natural numbers are presented in reflectica by the standard numerals, which are the strings $`0\q$, $`\S0\q$, $`\S\S0\q$ and so on. The constant symbols like $`1\q$ and $`2\q$ and expressions like $`1+2\q$ are reflectica expressions with expected values 
\[
\val(`1\q)=`\S0\q,\quad
\val(`2\q)=`\S\S0\q,
\quad\text{and}\quad
\val(`1+2\q)=`\S\S\S0\q. 
\]

Formulas are also expressions whose values could equal either $`\top\q$ (the true) or $`\bot\q$ (the false). For example, 
\[
\val(`1=1\q)=`\top\q
\quad\text{and}\quad
\val(`1=2\q)=`\bot\q.
\]
Note that, since every expression refers to its value and the expressions $`1=1\q$ and $`1=2\q$ have different values, we have
\[
\val(`(1=1)=(1=2)\q)=`\bot\q.
\]

Reflectica has means to refer to its own expressions. For example, $\val\big(`\sh(1+2)\q\big)=`1+2\q$ (the symbol $`{\sh}\q$ is a one-side quotation mark). Also reflectica allows you to work with quantifiers and set-theoretical notation (using formulas as sets, see Section~\ref{sect.conditional}). For example, 
\[
\val\big(`\exists{x}[({x}+1)=2]\q\big)=`\top\q,
\]
\[
\val\big(`1\in\sharp[({z}+1)=2)]\q\big)=`\top\q.
\]

To make the exposition of reflectica simpler and more uniform, we will use Polish (= prefix) notation everywhere. For example, we will write the expressions $1+{x}$ and $(1+{x})=2$ as $+(1,{x})$ and ${=}[{+}(1,{x}),2]$ respectively. Thus, expressions of reflectica have a form of ${s}\tau_1\!\ldots\tau_{n}$, where ${s}$ is an ${n}$-ary symbol, which we call the \textit{main symbol} of expression, and $\tau_1,\ldots,\tau_{n}$ are (sub)expressions, which we call the \textit{arguments} of expression ${s}\tau_1\!\ldots\tau_{n}$. The value of expression ${s}\tau_1\!\ldots\tau_{n}$, $\val({s}\tau_1\!\ldots\tau_{n})$, is determined by the rules for its main symbol ${s}$. These rules can be arbitrary, so the value of expression ${s}\tau_1\!\ldots\tau_{n}$ may depend not only on the objects $\val(\tau_1),\ldots,\val(\tau_{n})$, but also on the objects $\tau_1,\ldots,\tau_{n}$. If $\val({s}\tau_1\ldots\tau_{n})$ depends only on $\val(\tau_1),\ldots,\val(\tau_{n})$, we call ${s}$ an extensional symbol; otherwise we call ${s}$ an intensional symbol.

\section{The objects}
\noindent
The reflectica \textit{objects} are nonempty strings\footnote{A string is a finite sequence of symbols} in the reflectica \textit{alphabet}. We call the elements of this  alphabet reflectica \textit{symbols}. The values of reflectica expressions and the expressions themselves are objects. We identify each reflectica symbol with the string of length $1$ containing that symbol; therefore, symbols are also objects. The notions of a symbol and a string belong to the ontological base of reflectica, so we do not provide definitions for these notions. We assume that the reader understands the notions like \lq the first symbol of a string{\rq}, \lq the length of a string{\rq}  (that is, the number of symbols in the string), \lq the concatenation of two strings\rq , and so on. Thus, we do not define these notions either.

We use single quotes \lq`{\rq} and \lq'{\rq} in the usual way to denote reflectica strings. For example, `${+}{=}$' is the string of length 2 whose first symbol is the plus sign and the last symbol is the equality sign; we identify the first symbol the above string with the string `${+}$'.

To improve readability, we often use the punctuation marks comma \lq,{\rq} and brackets \lq(\rq, \lq)\rq, \lq[\rq, \lq]{\rq} inside the quotation marks \lq`{\rq} and \lq'\rq. Note that the above punctuation marks are not symbols of reflectica. But the symbols between quotes \lq`{\rq} and \lq'{\rq} --- with the exception of the five punctuation marks --- are always symbols of reflectica. For example, the string $`{=}({b},{2})\q$ contains 3 symbols and equals the string $`{=}{b}{2}\q$.

We denote the \textit{concatenation} of strings ${x}$ and ${y}$ by ${x}\cdot{y}$; here \lq${x}${\rq} and \lq${y}${\rq} are meta-variables (they may be thought of as variables of English) that we use to denote reflectica objects. For example, using this notation, we can assert that the string  `\(+\)'\(\cdot\)`\(=\)' equals the string `${+}{=}$'. We can write this statement  shorter: $`{+}\q\,{\cdot}\,`{=}\q\;=\;`{+}{=}\q$. Note that we used two different equality signs in the above statement: the middle one is a meta-language symbol that abbreviates English \lq equals{\rq} and the other two are reflectica symbols. We could add some accents or dots over symbols of reflectica (like in {\lq}$\dot{=}${\rq}) to distinguish between meta-language symbols and reflectica symbols, but we don't do this in order to relax notation. Instead we will put reflectica symbols between the quotes \lq`{\rq} and \lq'{\rq} in the contexts where ambiguity is possible.

\section{The canonical names of objects}\label{sect.canon.names}
\noindent
For every object ${x}$, reflectica has a special expression whose value equals ${x}$; we call that expression the canonical name of ${x}$. 
To construct canonical names, we use three symbols: `$\l$', `$\r$', and `$|$'; note that they are not punctuation marks. Using these symbols, we form an infinite sequence of pairing quotes: $`\l\q$ and $`\r\q$, $`|\l\q$ and $`\r|\q$, $`||\l\q$ and $`\r||\q$, \ldots\ . Then, to build the canonical name of a string ${x}$, we place around ${x}$ the first pair of the above quotes, whose closing quote is not a substring of ${x}$ (this method of quoting is similar to the unambiguous quotes of Boolos, see~\cite{Boolos1995-BOOQA}). Formally, we construct canonical names as follows:


Let
\[
{l}_0=`\l\q\quad\text{and}\quad{r}_0=`\r\q\,.
\]
For all natural ${n}\geq{0}$, let
\[
{l}_{{n}+1}=`|\q\cdot{l}_{n}\quad\text{and}\quad{r}_{{n}+1}={r}_{n}\cdot`|\q\,.
\]
For example, ${l}_1=`|\l\q$ and ${r}_3=`\r|||\q$. 

Let $\mathsf{m}({x})$ be the minimal natural number such that ${r}_{\mathsf{m}({x})}$ is not a substring of the string ${x}$. For example, $\mathsf{m}(`{+}\r|\q)=2$ and $\mathsf{m}(`\r{+}|\q)=1$.

We define \textit{the canonical name} of an object ${x}$, $\mathsf{name}({x})$, to be the string 
\[
{l}_{\mathsf{m}({x})}\cdot{x}\cdot{r}_{\mathsf{m}({x})}\ . 
\]
For example, 
\[
\mathsf{name}(`{+}{=}\q)=`\l{+}{=}\r\q,\quad
\mathsf{name}(`\l{+}{=}\r\q)=`|\l\,\l{+}{=}\r\,\r|\q,\quad
\mathsf{name}(`|\l\,\l{+}{=}\r\,\r|\q)=`||\l\,|\l\,\l{+}{=}\r\,\r|\,\r||\q.
\]
Note that $\name({x})\neq\name({y})$ whenever ${x}\neq{y}$. It follows that we can always recover an object from its canonical name (by removing the outer quotes). 

We say that an object ${c}$ is \textit{a canonical name} iff there is an object ${x}$ such that ${c}=\name({x})$. 
If ${c}$ is a canonical name, then no proper initial segment of ${c}$ is a canonical name.

\section{The expressions}
\noindent
The canonical names of objects are expressions of reflectica. We use standard Polish (= prefix) notation to build reflectica expressions. In accordance with this notation system, some reflectica symbols have the \textit{arity}, a natural number $\geq{0}$, which means the number of arguments. Note that the symbols `$\l$', `$\r$', and `$|$' do not have the arity. 

We define the notion of a reflectica \textit{expression} recursively as follows:
\begin{itemize}
  \item[E1.] If an object $\tau$ is a canonical name, then $\tau$ is an expression.
  \item[E2.] If the arity of a symbol ${s}$ equals ${n}$ and $\tau_1,\ldots,\tau_{n}$ are expressions, 
  then ${s}\cdot\tau_1\cdot\ldots\cdot\tau_{n}$ is an expression. (In particular, if a symbol ${s}$ has arity $0$, then ${s}$ is an expression.)
\end{itemize}
The arities of symbols $`{=}\q$ and $`\top\q$ are 2 and 0, respectively, so the following objects are expressions:
\[
`\top\q,\quad
`\l{+}{=}\r\q,\quad
`{=}\top\l{+}{=}\r\q.
\]

It is not hard to prove the unique readability lemma, see~\cite[Lemma II.4.3]{kunen2009foundations}, for the reflectica expressions. This lemma states that 
(1) if $\sigma$ is an expression, then no proper initial segment of $\sigma$ is an expression and 
(2) if the first symbol ${s}$ of an expression $\sigma$ has arity ${n}$, then $\sigma$ is not a canonical name and there are unique expressions $\tau_1,\ldots,\tau_{n}$ such that $\sigma$ equals ${s}\cdot\tau_1\cdot\ldots\cdot\tau_{n}$;
in this case we call the expressions $\tau_1,\ldots,\tau_{n}$ the \textit{arguments} of expression $\sigma$ and we call the symbol ${s}$ the \textit{main symbol} of expression $\sigma$. 

Note that it is possible to simplify the definition of <<reflectica expression>> by removing the base case (E1) and obtain a pure Polish notation system if the  alternative notation is used, see Section~\ref{sect.alt.expressions}.

\section{The value of an expression (the interpretation of reflectica)}
\label{val.of.term}
\noindent
The \textit{interpretation} of reflectica is a partial function $\val$ that maps expressions to their values, which are objects. We denote (in the meta-theory) the value of expression $\tau$ by $\val(\tau)$. For example, the interpretation maps the canonical name of an object ${x}$ to the object ${x}$, so $\val\big(\name({x})\big)={x}$ for every object ${x}$. We express the assertion $\val(\tau)={y}$ in three different ways: ${y}$ is the \textit{value} of $\tau$, $\tau$ is a \textit{name} for ${y}$, or $\tau$ \textit{refers to} ${y}$. An expression is \textit{defined} iff it has the value (that is, if the partial function $\val$ is defined on that expression).

Technically, $\val$ is a partial function from the set of  reflectica expressions to the universe of reflectica. We define the function $\val$ recursively. We call the statements that constitute the base of this recursion the \textit{axioms} of reflectica. We call the statements that constitute the steps of this recursion the \textit{rules} of reflectica. Reflectica has many axioms and rules, we will present them later. We say that reflectica is \textit{consistent} iff every its expression has no more than one value. So consistency of reflectica means that the recursive definition of the partial function $\val$ is correct. 

Here is an example of a reflectica axiom: $\val(`\top\q)=`\top\q$ (see Section~\ref{sect.prop.log.symb}). This axiom says that the expression $`\top\q$ is a name for the string $`\top\q$. 

Here is an example of a reflectica rule: For every object $\tau$, if $\tau$ is a formula and $\val(\tau)=`\top\q$, then $\val(`\neg\q\cdot\tau)=`\bot\q$ (this is an introduction rule for the symbol $`\neg\q$, see Section~\ref{sect.prop.log.symb}). 
We will write the rules in the traditional way: the premise of the rule is above the line and  the conclusion of the rule is under the line. In particular, we will write the above rule as follows:
\[
\frac{\tau\text{ is a formula and }\val(\tau)=`\top\q}{\val(`\neg\q\cdot\tau)=`\bot\q}\ .
\]

Above we say that the canonical name of an object ${x}$ refers to the object ${x}$. But the values of expressions should be defined by axioms and rules alone, so we add to reflectica the following axiom scheme:
\[
\val\big(\name({x})\big)={x}\quad\text{for every object }{x}.
\]

As an example of using this recursive definition of $\val$, we can derive from the above axioms and rule the equalities $\val(`\top\q)=`\top\q=\val(`\l\top\r\q)$ and $\val(`\neg\top\q)=`\bot\q$  (the string $`{\top}\q$ is a formula, see Sections~\ref{sect.formulas} and~\ref{sect.prop.log.symb}).

Informally, when $\val({s}\,{\cdot}\,\tau_1{\cdot}\ldots{\cdot}\,\tau_{n})$ depends only on $\val(\tau_1),\ldots,\val(\tau_{n})$, we call ${s}$ an extensional symbol; otherwise we call ${s}$ an intensional symbol. 

\section{The meaning of a symbol (the semantics of reflectica)}
\label{sect.meaning}
\noindent
In an expression of the form ${s}\,{\cdot}\,\tau_1{\cdot}\ldots{\cdot}\,\tau_{n}$ --- where ${s}$ is an ${n}$-ary symbol --- we treat the main symbol ${s}$ as an operation on the arguments $\tau_1,\ldots,\tau_{n}$, with $\val({s}\,{\cdot}\,\tau_1{\cdot}\ldots{\cdot}\,\tau_{n})$ being the result of such an operation.
From this point of view, the \textit{meaning} of the symbol ${s}$ is the operation it performs. The result of such an operation is determined by the axioms and rules for the symbol ${s}$. Therefore, in reflectica, the meaning of a symbol is determined by the axioms and rules for that symbol. The symbols of reflectica correspond to the words of natural language. Thus, in a sense, reflectica reconciles the perspective that the meaning of a word is determined by the values (of expressions in which that word is the main word) with the perspective that the meaning of a word is determined by the rules of its use.

To be a self-reflecting language, the reflectica language must have expressive means to describe its own semantics, which consists in the meanings of reflectica symbols. In Section~\ref{sect.dagger} we show that reflectica can express assertions about its own expressions and their values. In Section~\ref{sect.rules} we show that reflectica can express the statements made by its axioms and has expressive means to declare its own rules. Thus, from both of the above-mentioned perspectives on what the meanings of words/symbols are, the reflectica language allows us to describe its own semantics.

\section{The formulas}\label{sect.formulas}
\noindent
А logical formula is an expression that expresses some statement, so its intended value is either true or false. In extensional (i.e., not intensional) systems~\cite{sep-logic-intensional}, a formula whose main symbol is an $n$-ary predicate symbol is an assertion about $n$ objects that are the values of the arguments. Some formulas of extensional systems are not such statements --- for example, a formula whose main symbol is a quantifier (which is a $2$-ary symbol) is not a statement about two objects that are the values of the arguments, so quantifiers are not predicate symbols in such systems.

In reflectica, the arguments of an expression are also objects, so the main symbol (of a formula) can act as an assertion about objects that include not only the values of the arguments, but also the arguments themselves. It is quite natural to call such symbols (intensional) predicate symbols. In particular, we may consider logical symbols such as $`\&\q$, $`{\to}\q$, or $`\forall\q$ as (intensional) predicate symbols of reflectica. For example, the reflectica symbol $`\forall\q$ says that its first argument is a variable and if we replace the free occurrences of that variable in the second argument with an arbitrary expression, the result will be a formula with value the true. Therefore, we may simplify the standard notion of a formula as follows:

\smallskip
We say that a reflectica expression $\varphi$ is a reflectica \textit{formula} iff  the main symbol of $\varphi$ is a predicate symbol. Below we introduce predicate symbols, among which are the standard logical symbols.
We say that a reflectica formula $\varphi$ is a reflectica \textit{theorem} iff $\val(\varphi)=`\top\q$.
We call the symbol $`\top\q$ 
the \textit{true},
it has arity 0. So a theorem is a formula whose value equals the true. 

\smallskip
Informally\footnote{In order to give a formal definition here, it is necessary to specify the meta-theory we use to construct reflectica.}, when we say that a reflectica formula $\varphi$ \textit{represents} a meta-language statement ${S}$, we mean that $\varphi$ is a reflectica theorem whenever the statement ${S}$ holds; when we say that a formula $\varphi$ \textit{declares} a statement ${S}$, we mean that if $\varphi$ is a theorem, then the statement ${S}$ holds;
and when we say that a formula $\varphi$ \textit{expresses} a statement ${S}$, we mean that $\varphi$ represents and declares the statement ${S}$. 

For example, if an 1-ary predicate symbol $`{\bullet}\q$ has the (introduction) rule 
\[
\frac{\tau\text{\ is an expression and\ } \val(\tau)=`\forall\q}{\val(`{\bullet}\q\cdot\tau)=`\top\q},
\]
then a formula of the form $`{\bullet}\q\cdot\tau$ represents the statement that the value of its argument $\tau$ equals $`{\forall}\q$.

If an 1-ary predicate symbol $`{\circ}\q$ has the (elimination) rule 
\[
\frac{\tau\text{\ is an expression and\ } \val(`{\circ}\q\cdot\tau)=`\top\q}{\val(\tau)=`\forall\q},
\]
then a formula of the form $`{\circ}\q\cdot\tau$ declares the statement that the value of its argument $\tau$ equals $`{\forall}\q$.

And if an 1-ary predicate symbol $`{\star}\q$ has the rules
\[
\frac{\tau\text{\ is an expression and\ } \val(\tau)=`\forall\q}{\val(`{\star}\q\cdot\tau)=`\top\q}
\quad\text{and}\quad
\frac{\tau\text{\ is an expression and\ } \val(`{\star}\q\cdot\tau)=`\top\q}{\val(\tau)=`\forall\q},
\]
then a formula of the form $`{\star}\q\cdot\tau$ expresses the statement that the value of its argument $\tau$ equals $`{\forall}\q$.

Informally, when we say that reflectica \textit{knows} a meta-language statement ${S}$, we mean that there exists a reflectica theorem that expresses the statement ${S}$. 
For example, if we add to reflectica the above symbol $`{\star}\q$ with its two rules, then the formula $`{\star}\l{\forall}\r\q$ will be a theorem, so reflectica will know that $\val(`\l{\forall}\r\q)=`{\forall}\q$ (in fact, reflectica knows that $\val(`\l{\forall}\r\q)=`{\forall}\q$ even without the symbol $`{\star}\q$ because this assertion is an axiom and reflectica knows the statements made by its axioms, see Section~\ref{sect.rules}).

\section{The symbols of propositional logic}\label{sect.prop.log.symb}
\noindent
The symbols of propositional logic are predicate symbols of reflectica; among them are the \textit{conjunction} $`\&\q$ and the \textit{disjunction} $`{\vee}\q$ of arity 2; the \textit{negation} $`\neg\q$ of arity 1; the \textit{true} $`\top\q$ and the \textit{false} $`\bot\q$ of arity 0. Their axioms and rules correspond the truth-tables of Kleene's 3-valued strong logic of indeterminacy $K_3$~\cite{sep-truth-values}. Note that not every formula of reflectica has a value.

The symbols $`\top\q$ and $`\bot\q$ have the following axioms:
\[
\val(`\top\q)=`\top\q
\qquad\text{ and }\qquad
\val(`\bot\q)=`\bot\q.
\]
The symbols $`\&\q$, $`{\vee}\q$, and $`\neg\q$ have the following introduction rules:
\[
\frac{\phi\text{\ is a formula and\ } \val(\phi)=`\top\q}{\val(`\neg\q\cdot\phi)=`\bot\q}\ ,
\qquad
\frac{\phi\text{\ is a formula and\ } \val(\phi)=`\bot\q}{\val(`\neg\q\cdot\phi)=`\top\q}\ ,
\]
\[
\frac{\phi,\psi\text{\ are formulas  and\ }\val(\phi)=\val(\psi)=`\top\q}{\val(`\&\q\cdot\phi\cdot\psi)=`\top\q}\ ,
\]
\[
\frac{\phi,\psi\text{\ are formulas  and\ }\val(\phi)=`\bot\q}{\val(`\&\q\cdot\phi\cdot\psi)=`\bot\q}\ ,
\qquad
\frac{\phi,\psi\text{\ are formulas  and\ }\val(\psi)=`\bot\q}{\val(`\&\q\cdot\phi\cdot\psi)=`\bot\q}\ ,
\]
\[
\frac{\phi,\psi\text{\ are formulas  and\ }\val(\phi)=\val(\psi)=`\bot\q}{\val(`{\vee}\q\cdot\phi\cdot\psi)=`\bot\q}\ ,
\]
\[
\frac{\phi,\psi\text{\ are formulas  and\ }\val(\phi)=`\top\q}{\val(`{\vee}\q\cdot\phi\cdot\psi)=`\top\q}\ ,
\qquad
\frac{\phi,\psi\text{\ are formulas  and\ }\val(\psi)=`\top\q}{\val(`{\vee}\q\cdot\phi\cdot\psi)=`\top\q}\ .
\]

\section{The equality symbol `${=}$'}\label{sect.=}
\noindent
The symbol $`{=}\q$, a predicate symbol of arity 2, allows reflectica to make assertions about equality of objects: a formula of the form $`{=}\q\cdot\tau\cdot\sigma$ says that the expressions $\tau$ and $\sigma$ have values and these values are equal. Here are introduction rules for $`{=}\q$: 
\begin{prooftree}
\def\fCenter{\mbox{\ $\to$\ }}
\AxiomC{$\tau,\sigma$ are defined expressions and $\val(\tau)=\val(\sigma)$}
\RightLabel{\ ,}
\UnaryInfC{$\val(`{=}\q\cdot\tau\cdot\sigma)=`\top\q$}
\end{prooftree}
\begin{prooftree}
\def\fCenter{\mbox{\ $\to$\ }}
\AxiomC{$\tau,\sigma$ are defined expressions and $\val(\tau)\neq\val(\sigma)$}
\RightLabel{\ .}
\UnaryInfC{$\val(`{=}\q\cdot\tau\cdot\sigma)=`\bot\q$}
\end{prooftree}
For example, using these rules, it can be deduced that $\val(`{=}\l\forall\r\l\forall\r\q)=`\top\q$ and 
$\val(`{=}\l\forall\r\l\exists\r\q)=`\bot\q$. So $`{=}\l\forall\r\l\forall\r\q$ and the negation of $`{=}\l\forall\r\l\exists\r\q$ are reflectica theorems.

Here are elimination rules for equality:
\begin{prooftree}
\def\fCenter{\mbox{\ $\to$\ }}
\AxiomC{$\tau,\sigma$ are expressions, $\val(\tau)={x}$, and $\val(`{=}\q\cdot\tau\cdot\sigma)=`\top\q$}
\RightLabel{\ ,}
\UnaryInfC{$\val(\sigma)={x}$}
\end{prooftree}
\begin{prooftree}
\def\fCenter{\mbox{\ $\to$\ }}
\AxiomC{$\tau,\sigma$ are expressions, $\val(\sigma)={x}$, and $\val(`{=}\q\cdot\tau\cdot\sigma)=`\top\q$}
\RightLabel{\ .}
\UnaryInfC{$\val(\tau)={x}$}
\end{prooftree}

\section{One-side quotation marks `$\sh$' and `$\natural$'}\label{sect.quotes}
\noindent
The (intensional) symbol $`\sh\q$ of arity 1 provides an easy way to refer to expressions: if $\tau$ is an expression, then the expression $`\sh\q\cdot\tau$ is a name for $\tau$. 
The axiom scheme for $`\sh\q$ is simple:
\[
\val(`\sh\q\cdot\tau)=\tau\quad\text{for every expression }{\tau}.
\]
For example, the value of the expression $`{\sharp}\l\neg\r\q$ equals $`\l\neg\r\q$.

The (intensional) symbol $`\na\q$ of arity 1 is another one-side quote for expressions, but, unlike $`\sh\q$, $`\na\q$ does not bind variables inside its argument, see the next section. An expression of the form $`\natural\q\cdot\tau$ refers to the expression $\tau$, which is postulated by the axiom scheme
\[
\val(`\natural\q\cdot\tau)=\tau\quad\text{for every expression }{\tau}.
\]

\section{Free and bound variables}
\noindent
Some symbols of reflectica are \textit{variables}, they have arity 0, so they are expressions. In particular, the letters $`{a}\q,\ldots,`{z}\q$ and $`\alpha\q,\ldots,`\omega\q$  of Latin and Greek alphabets are reflectica variables. 

As in classical predicate logic, every occurrence of a variable in an expression is either \textit{free} or \textit{bound}; precise definition of these notions can be found in any textbook on predicate logic. Briefly, it can be formulated as follows: an occurrence of a variable ${v}$ in an expression is bound iff it lies inside the scope of a quantifier $`\forall\q$ or $`\exists\q$ acting on (i.e., followed by) ${v}$, see~\cite[Definition II.5.5]{kunen2009foundations}.  
We extend the classical definition by adding that quoting --- with the exception of the one-side quote $`\na\q$ --- binds all variables inside the quotation:
\begin{itemize}
  \item[V1.] Every occurrence of a variable inside a canonical name is bound. 
  \item[V2.] The one-side quote $`{\sharp}\q$ binds all variables in its argument.
\end{itemize}

\section{The quantifiers `$\exists$' and `$\forall$'}\label{sect.quantifiers}
\noindent
The \textit{existential quantifier} $`\exists\q$ and the \textit{universal quantifier} $`\forall\q$ are (intensional) predicate symbols of arity 2. 
A formula of the form $`\exists\q\cdot{v}\cdot\phi$ represents the statement that its first argument is a variable and there exists an expression $\tau$ such that if we replace all free occurrences of that variable in the second argument with $\tau$, the result will be a theorem.
A formula of the form $`\forall\q\cdot{v}\cdot\phi$ declares the statement that its first argument is a variable and if we replace all free occurrences of that variable in the second argument with an arbitrary expression $\tau$, the result will always be a theorem.

To write the rules for quantifiers, we need the following notation. For a variable ${v}$ and  expressions $\phi$ and $\tau$, we denote by $\mathsf{Sub}(\varphi,{v}\leadsto\tau)$ the expression which results from $\varphi$ by replacing all free occurrences of variable ${v}$ by expression $\tau$. Note that this is a meta-theoretic notation and \lq${v}$\rq, \lq$\phi$\rq, \lq$\tau${\rq} are meta-variables that denote objects of reflectica. 
For example, 
\[
\mathsf{Sub}\big(`\&[{=}({x},{y}),\forall{y}{=}({x},{y})]\q,`{y}\q\leadsto`\l{+}\r\q\big)\quad\text{equals the string}\quad`\&[{=}({x},\l{+}\r),\forall{y}{=}({x},{y})]\q.
\]
(In more common notation the formula $`\&[{=}({x},{y}),\forall{y}{=}({x},{y})]\q$ would be written as ${x}={y}\ \&\ \forall{y}({x}={y})$.)

Here is the introduction rule for the existential quantifier:
\begin{prooftree}
\def\fCenter{\mbox{\ $\to$\ }}
\AxiomC{${v}$ is a variable, $\varphi$ and $\tau$ are expressions, 
and $\mathsf{Sub}(\varphi,{v}\leadsto\tau)$ is a theorem}
\RightLabel{\ .}
\UnaryInfC{$\val(`\exists\q\cdot{v}\cdot\varphi)=`\top\q$}
\end{prooftree}

Here is the elimination rule for the universal quantifier:
\begin{prooftree}
\def\fCenter{\mbox{\ $\to$\ }}
\AxiomC{${v}$ is a variable, $\varphi$ and $\tau$ are expressions, and $\val(`\forall\q\cdot{v}\cdot\varphi)=`\top\q$}
\RightLabel{\ .}
\UnaryInfC{$\val\big(\mathsf{Sub}(\varphi,{v}\leadsto\tau)\big)=`\top\q$}
\end{prooftree}

As an example, notice that $`\exists{x}{x}\q$ is a reflectica theorem. As a more usual example, consider the following reflectica formula, which states that for every natural number, its successor is not 0:
\[
`\:\forall{n}\,{\to}\,\big(\mathsf{\underline{Number}}({n}),\neg{=}\,[\S({n}),0]\big)\:\q
\]
--- in this formula, the 1-ary extensional predicate symbol $`\mathsf{\underline{Number}}\q$ says that its argument refers to a natural number (see Section~\ref{sect.numbers}) and $`{\to}\q$ is the conditional from Section~\ref{sect.conditional}.

\section{The concatenation of strings}
\noindent
The symbol $`{\cdot}\q$ of arity 2 acts as the (extensional) concatenation operation. An expression of the form $`{\cdot}\q\cdot\tau\cdot\sigma$ refers to the concatenation of the values of its arguments $\tau$ and $\sigma$. The introduction rule for $`{\cdot}\q$ is straightforward: 
\begin{prooftree}
\def\fCenter{\mbox{\ $\to$\ }}
\AxiomC{$\tau,\sigma$ are defined expressions and $\val(\tau)\cdot\val(\sigma)={x}$}
\RightLabel{\ .}
\UnaryInfC{$\val(`{\cdot}\q\cdot\tau\cdot\sigma)={x}$}
\end{prooftree}
For example, $\val(`\,{\cdot}\,\l{+}\r\l{=}\r\q)=`{+}{=}\q$.

\section{The natural numbers}
\label{sect.numbers}
\noindent
The symbols $`0\q$ and $`\S\q$ have arities 0 and 1, respectively. Using these symbols, we can write the expressions $`0\q$, $`\S0\q$, $`\S\S0\q$ and so on, which are the standard numerals; we call them the \textit{natural numbers} of reflectica.

The expression $`0\q$ refers to itself, this is postulated by the axiom $\val(`0\q)=`0\q$. 
If the value of an expression $\tau$ is a natural number, then the expression $`\S\q\cdot\tau$ refers to the successor of that natural number. Accordingly, we have the  introduction rule
\[
\frac{\tau\text{ is an expression and }\val(\tau)={x}}{\val(`\S\q\cdot\tau)=`\S\q\cdot{x}}\ .
\]
It follows that $\val(`0\q)=`0\q$, $\val(`\S0\q)=`\S0\q$, $\val(`\S\S0\q)=`\S\S0\q$ and so on; that is, numerals refer to the corresponding natural numbers, which is in accordance with common practice.

\section{The symbol `\hspace{0.5pt}\textmd{\textdagger}\hspace{0.5pt}' that allows reflectica to make assertions about its own interpretation}\label{sect.dagger}
\noindent
The interpretation of reflectica is the partial function $\val$, which assigns to expressions their values. To be a self-reflecting language, the language of reflectica must have expressive means to describe its own interpretation; that is, to make assertions about expressions and their values. Reflectica can mention expressions using their names. 
To mention the value of an expression, reflectica has a symbol $`\dagger\q$, the dagger, such that the operation $`\dagger\q$ acts exactly like the function $\val$. 
This symbol an extensional 1-ary symbol and an expression of the form $`\dagger\q\cdot\tau$ refers to the value of the expression that is the value of the argument $\tau$. 

For example, the expression $`\dagger{\cdot}(\l\neg\r,\l\top\r)\q$ has the argument 
$`{\cdot}(\l\neg\r,\l\top\r)\q$; this arguments value is the expression $`\neg\top\q$, whose value equals the string $`\bot\q$; therefore the expression $`\dagger{\cdot}(\l\neg\r,\l\top\r)\q$ refers to the string $`\bot\q$. 
So the true meta-language statement
\[
\val(`\neg\q{\cdot}`\top\q)=`\bot\q
\] 
can be expressed in reflectica by the formula (actually, a theorem)
\[
`\:{=}\,\big(\dagger{\cdot}(\l\neg\r,\l\top\r),\l\bot\r\big)\:\q.
\] 
In other words, the dagger acts as the operation \lq the value of the value of the argument{\rq}. Accordingly, reflectica has the following introduction rule for $`\dagger\q$: 
\[
\frac{\tau,\sigma\text{\ are expressions,\ } \val(\tau)=\sigma,\ \val(\sigma)={x}}
{\val(`\dagger\q\cdot\tau)={x}}\ .
\]
Note that there is no need for a symbol that acts as the operation \lq the value of the argument{\rq} since that value is already referred to by the argument itself (but we may introduce such a symbol if we want).

Here is the elimination rule for $`\dagger\q$:
\[
\frac{\tau,\sigma\text{\ are expressions,\ } \val(\tau)=\sigma,\ \val(`\dagger\q\cdot\tau)={x}}
{\val(\sigma)={x}}\ .
\]

\section{The symbol `$\rightY$' that allows reflectica to declare its rules}
\label{sect.rules}
\noindent
Reflectica can express assertions made by its axioms. For example, the formula $`{=}(\dagger\l\bot\r,\l\bot\r)\q$ expresses the statement $\val(`\bot\q)=`\bot\q$ (check it by using the introduction and elimination rules for $`{=}\q$ and $`{\dagger}\q$), which is an axiom of reflectica. It is easy to check that the above formula is a theorem, so reflectica {\lq}knows{\rq} (see the end of Section~\ref{sect.formulas}) the statement made by the  above axiom.\footnote{Note that the same statement can be expressed by a simpler formula $`{=}(\bot,\l\bot\r)\q$.} Similarly, reflectica knows all statements that are asserted by its axioms of the form $\val(\tau)={x}$.

To make assertions about its rules, reflectica has the symbol $`{\rightY}\q$, an (intensional) predicate symbol of arity 2. A formula of the form $`{\rightY}\q\cdot\varphi\cdot\psi$ declares that reflectica obeys the rule whose premise is expressed by the formula $\varphi$ and whose conclusion is expressed by the formula $\psi$. 
Here is the elimination rule for $`{\rightY}\q$:
\[
\frac{\varphi,\psi\text{\ are formulas,\ }\val(`{\rightY}\q{\cdot}\varphi{\cdot}\psi)\,{=}\,`\top\q, \text{\ and\ }\val(\varphi)\,{=}\,`\top\q}
{\val(\psi)=`\top\q}\,.
\]

For example, using $`{\rightY}\q$ and the 1-ary (extensional) predicate symbol $`\mathsf{\underline{Form}}\q$ (such that the formula $`\mathsf{\underline{Form}}\q\cdot\tau$ expresses the statement that $\val(\tau)$ is a formula), we can declare the rule
\[
\frac{\varphi\text{\ is a formula and\ } \val(\varphi)=`\top\q}{\val(`\neg\q\cdot\varphi)=`\bot\q}
\]
by the formula 
\[
\begin{split}
`\ \forall{\varphi}\ {\rightY}\ &\&\,[\hspace{0.5pt}\mathsf{\underline{Form}}(\na\hspace{0.5pt}\varphi),{=}\,(\varphi,\top)]\\
 &{=}\,(\neg\varphi,\bot)\ \q
\end{split}
\]
(recall that $`\na\q$ is the one-side quote that does not bind variables, $`\varphi\q$ is a variable, $\val(`\top\q)=`\top\q$, and $\val(`\bot\q)=`\bot\q$). 
If we want to declare the above rule by using reflectica language, it is enough to add the axiom 
\[
\val(`\forall\varphi{\rightY}\&\mathsf{\underline{Form}}\na\varphi{=}\varphi\top{=}\neg\varphi\bot\q)=`\top\q.
\]

In a similar way, we can declare a new rule by adding a new axiom to reflectica.  
For example, we could declare the rule
\[
\frac{\tau\text{\ is an expression and\ } \val(\tau)=`\forall\q}{\val(`\neg\q\cdot\tau)=`\exists\q}
\]
by adding the axiom 
\[
\val\big(`\,\forall\tau{\rightY}[{=}(\tau,\l\forall\r),{=}(\neg\tau,\l\exists\r)]\:\q\big)=`\top\q\ .
\]
It is important to note that reflectica will obey the above rule after adding such an axiom, even if we do not add the above rule to reflectica (check it by using the elimination rules for $`{\forall}\q$ and $`{\rightY}\q$ and introduction and elimination rules for $`{=}\q$). Thus, if we want to extend reflectica (for example, by adding a new symbol and rules for it), there is no need to add new rules to reflectica. Instead, we can add new axioms (in the form $\val(\phi)=`\top\q$) that declare the rules we need.

In a similar way we can declare (and add to reflectica if we want) an axiom scheme with just one reflectica formula. For example, the axiom scheme 
\[
\val(`\sh\q\cdot\tau)=\tau\quad\text{for every expression }{\tau}
\]
from Section~\ref{sect.quotes} can be declared by the reflectica formula
\[
`\,\forall\tau\,{\to}\,\big(\mathsf{\underline{Expr}}(\tau),{=}\,[\dagger{\cdot}(\l\sh\r,\tau),\tau]\big)\,\q\ ,
\]
where the 1-ary extensional predicate sybmol $`\mathsf{\underline{Expr}}\q$ says that the value of its argument is an expression and $`{\to}\q$ is the conditional from Section~\ref{sect.conditional}. 

\smallskip

If we want all reflectica rules to be declared by its axioms, we should write formulas declaring that rules and add to reflectica new axioms (in the form $\val(\phi)=`\top\q$) which say that these formulas are theorems. Note that, as we mentioned in the beginning of this section, reflectica will know the statements asserted by these new axioms as soon as we add them. 

To implement this plan, it is convenient first to add symbols that act as the operations like {\lq}the arity of a symbol{\rq}, {\lq}the main symbol of an expression{\rq}, and so on and predicate symbols that express the properties of {\lq}being an expression{\rq}, {\lq}an occurrence of a variable in an expression being free{\rq}, etc. The axioms and rules for such symbols are quite straightforward. 

\section{Conditional sentences and Curry's paradox}\label{sect.conditional}
\noindent
To express conditional sentences --- that is, sentences of the form ``if A, then B'' --- 
classical logic uses implication, which can be defined through negation and disjunction. Following this pattern, we add a 2-ary predicate symbol $`{\Rightarrow}\q$, the \textit{material implication}, with an introduction rule
\[
\frac{\phi,\psi\text{\ are formulas  and\ }\val(`{\vee}\neg\q\cdot\phi\cdot\psi)={x}}
{\val(`{\Rightarrow}\q\cdot\phi\cdot\psi)={x}}\ .
\]
The problem is that an expression of the form 
$`{\Rightarrow}\q\cdot\phi\cdot\psi$ is defined only if at least one of its arguments is defined. Therefore, this symbol is not a good choice for conditional sentences.

To write conditional sentences, we introduce another 2-ary predicate symbol, $`{\to}\q$, the \textit{conditional}. We declare rules for it by adding axioms (so we do not add new rules, see the previous section) that informally can be written as follows:
\begin{gather*}
(\varphi\rightY\psi)\ \rightY\ (\varphi\to\psi),
\\    
(\varphi\Rightarrow\psi)\ \rightY\ (\varphi\to\psi),
\\
(\varphi\to\psi)\,\&\,(\psi\to\xi)\ \rightY\ (\varphi\to\xi),
\\
(\varphi\to\psi)\,\&\,(\phi=\top)\ \rightY\ (\psi=\top).
\end{gather*}
Formally, the first of these axioms could be added to reflectica in the form 
\[\val\big(`\forall\varphi\forall\psi\rightY[\rightY(\varphi,\psi),{\to}(\varphi,\psi)]\q\big)=`\top\q.
\]

We should be careful here because if we add too many or too strong axioms, we could face Curry's paradox. So we should work in the style of substructural logics~\cite{sep-logic-substructural} here. For example, if we add the modus ponens rule in the form of the axiom (informally written as)
\[
\big[(\varphi\to[\psi\to\xi])\,\&\,
(\varphi\to\psi)\big]\ \rightY\ (\varphi\to\xi),
\]
and also add a 2-ary predicate symbol $`{\in}\q$ such that a formula of the form $`{\in}\q\cdot\tau\cdot\sigma$ says that $\mathsf{Sub}(\val(\sigma),{v}\leadsto\tau)$ is a theorem, where ${v}$ is the free variable of the formula $\val(\sigma)$, then we could get Curry's paradox, where Curry sentence~\cite{sep-curry-paradox} is the formula
\[
([\varphi\in\natural\varphi]\to\bot)\,\in\,
\natural([\varphi\in\natural\varphi]\to\bot).
\]
But it looks safe to add the following weakened form of modus ponens:
\[
\big[(\varphi\to[\psi\to\xi])\,\&\,
(\varphi\to\psi)\,\&\,({\downarrow}\varphi)\big]\ \rightY\ (\varphi\to\xi),
\]
where $`{\downarrow}\q$ is an (intensional) 1-ary predicate symbol such that a formula of the form $`{\downarrow}\q\cdot\tau$ says that its argument $\tau$ is defined.

Using the conditional symbol, we may add to reflectica an axiom saying that every reflectica formula is equivalent to the truthness of its value, which can be informally written as
\[
\forall{\varphi}\:\big[\hspace{0.5pt}\mathsf{\underline{Form}}(\na\hspace{0.5pt}\varphi)\to\big([(\varphi=\!\top)\,{\to}\,\varphi]\,\&\,[\varphi\,{\to}\,(\varphi=\!\top)]\big)\big]\ .
\]

\section{Intensional context in reflectica}\label{sect.context}
\noindent
Reflectica allows us to express statements containing an intensional context (see footnote on page \pageref{footnote context}); this ability makes reflectica an intensional system~\cite{sep-logic-intensional}.  For example, we may introduce a 2-ary predicate symbol $`\mathsf{\underline{Believes}}\q$ that expresses the corresponding propositional attitude~\cite{sep-prop-attitude-reports} so that a formula of the form $`\mathsf{\underline{Believes}}\q\cdot\tau\cdot\varphi$ says that the object $\val(\tau)$ believes the statement expressed by the formula $\varphi$. Note that such a symbol $`\mathsf{\underline{Believes}}\q$ is extensional in the first argument and intensional in the second. To improve readability, we write most of reflectica expressions in this section informally.

Using such a symbol, we can say that every statement believed by the object named by the symbol $`\mathsf{\underline{Leo}}\q$ is a theorem:
\[
\forall\varphi\,[\mathsf{\underline{Form}}(\natural\varphi)\,\&\,\mathsf{\underline{Believes}}(\mathsf{\underline{Leo}},\varphi)\;\to\;\varphi].
\]
Or we can say that Leo believes every  reflectica rule:
\[
\forall\varphi\forall\psi\,[(\phi\,{\rightY}\,\psi)\ \to\ \mathsf{\underline{Believes}}(\mathsf{\underline{Leo}},\phi\,{\rightY}\,\psi)].
\]

We may add to reflectica a rule declaring that everyone believes every statement that follows by a rule he believes from a statement he believes:
\[
\forall{x}\forall\varphi\forall\psi\,\big[\mathsf{\underline{Believes}}({x},\varphi)\,\&\,
\mathsf{\underline{Believes}}({x},\varphi\rightY\psi)\ \rightY\ 
\mathsf{\underline{Believes}}({x},\psi)\big].
\]

As another example of an intensional context, reflectica allows us to correctly distinguish between \textit{de re} and \textit{de dicto} meanings of a sentence, see the Supplement to~\cite{sep-prop-attitude-reports}. For example, the sentence \textit{Leo believes that some number is prime} can mean either
\[
\mathsf{\underline{Believes}}\big(\mathsf{\underline{Leo}},\exists{x}[\mathsf{\underline{Number}}({x})\&\mathsf{\underline{Prime}}({x})]\big)
\]
or
\[
\exists{x}\big(\mathsf{\underline{Number}}({x})\:\&\:
\mathsf{\underline{Believes}}[\mathsf{\underline{Leo}},\mathsf{\underline{Prime}}({x})]\big).
\]
Note that, since the symbol $`\mathsf{\underline{Believes}}\q$ is intensional in the second argument, the latter formula involves quantifying into an intensional context, which Quine thought is incoherent~\cite{quine1956quantifiers} (but reflectica allows to express such things coherently). 

As another example of quantifying into an intensional context, reflectica allows us to correctly quantify into a quotation by using the one-side quotation symbol $`{\natural}\q$ (because this symbol does not bind variables). For example, consider the formula saying that there exists an expression ${x}$ such that the formula which says this expression refers to a prime number, has true value:
\[
`\,\exists{x}{=}\big(\dagger\natural\&[\mathsf{\underline{Number}}({x}),\mathsf{\underline{Prime}}({x})],\!\top\big)\:\q.
\]
In this formula, the variable $`{x}\q$ is inside the quote expression $`\natural\&[\mathsf{\underline{Number}}({x}),\mathsf{\underline{Prime}}({x})]\q$ and at the same time it is not trapped by the quotation mark $`{\natural}\q$. Some philosophers say that \lq it is not possible to quantify into quotation\rq, see BQ2 in~Section 2 of \cite{sep-quotation}.

\section{Further development}
\noindent
We need to add symbols, axioms and rules to reflectica to make it a full-fledged theory; for example, we may add an introduction rule for the universal quantifier in the form
\begin{prooftree}
\def\fCenter{\mbox{\ $\to$\ }}
\AxiomC{${v}$ is a variable, $\varphi$ is an expressions, and for every expression $\tau$, 
$\mathsf{Sub}(\varphi,{v}\leadsto\tau)$ is a theorem}
\RightLabel{\ .}
\UnaryInfC{$\val(`\forall\q\cdot{v}\cdot\varphi)=`\top\q$}
\end{prooftree}
This can be implemented in many various ways, leading to various theories. Such theories will be explored in future articles.

\section{Appendix: an alternative pure Polish notation}
\label{sect.alt.expressions}
\noindent

We could simplify the definition of the notion <<reflectica expression>> and use pure Polish notation as follows. Instead of using canonical names for objects, we could use <<canonical names for symbols>>: for each reflectica symbol, we can add the paired symbol --- its \textit{check name} ---  that we can draw in the same way as the original symbol, but with a check on top. For example, the symbol $`\check{\forall}\q$ is the check name for the symbol $`{\forall}\q$. We will not add check names for these new symbols (in other words, we do not add symbols with more then one check on top). 

Let us denote (in the meta-theory) the check name for the symbol ${s}$ by $\mathsf{chn}({s})$. Each check name is a 0-ary symbol and $\val\big(\mathsf{chn}({s})\big)={s}$ for every symbol ${s}$ that is not a check name.

In such a setting, every reflectica symbol gets a name: if a symbol ${s}$ is not a check name, then $\mathsf{chn}({s})$ refers to it; if ${s}$ is a check name, then the expression $`\sh\q\cdot{s}$ is a name for ${s}$ (for example, $`\sh\check{\forall}\q$ is a name for $`\check{\forall}\q$). We can not use $`\sh\q\cdot{t}$ as a name for an arbitrary symbol ${t}$, because $`\sh\q\cdot{t}$ is an expression only in case ${t}$ is a 0-ary symbol.

Once we have a name for each symbol, we can give a name for each string by simply using concatenation and
symbol names. For example, the object $`\,{+}\check{=}{=}\,\q$ gets a name $`\,{\cdot}\,[{\cdot}\,(\check{+},\sh\check{=}),\check{=}]\,\q$. In such a setting, we can do without canonical names and get a pure Polish notation with the following simplified (recursive) definition of reflectica \textit{expression}:
\begin{itemize}
  \item[E.] If the arity of a symbol ${s}$ equals ${n}$ and $\tau_1,\ldots,\tau_{n}$ are expressions, 
  then ${s}\cdot\tau_1\cdot\ldots\cdot\tau_{n}$ is an expression. (In particular, if a symbol ${s}$ has arity $0$, then ${s}$ is an expression.)
\end{itemize}

We do not use this alternative notation in this paper because its worse readability: compare the canonical name 
\[
\l{+}\check{=}{=}{+}\check{=}{=}\r
\] 
for the string
\[
{+}\check{=}{=}{+}\check{=}{=}
\]
with its name 
\[
{\cdot}\,{\cdot}\,{\cdot}\,{\cdot}\,{\cdot}\,\check{+}\sh\check{=}\check{+}\check{+}\sh\check{=}\check{+}
\] 
in the alternative notation.


\bibliographystyle{plain} 
\bibliography{Reflectica.bib}

\begin{thebibliography}{10}

\bibitem{Boolos1995-BOOQA}
George Boolos.
\newblock Quotational ambiguity.
\newblock In Paolo Leonardi and Marco Santambrogio, editors, {\em On Quine: New Essays}, pages 283--296. Cambridge University Press, 1995.

\bibitem{sep-quotation}
Herman Cappelen, Ernest Lepore, and Matthew McKeever.
\newblock {Quotation}.
\newblock In Edward~N. Zalta and Uri Nodelman, editors, {\em The {Stanford} Encyclopedia of Philosophy}. Metaphysics Research Lab, Stanford University, {S}ummer 2023 edition, 2023.

\bibitem{sep-logic-intensional}
Melvin Fitting.
\newblock {Intensional Logic}.
\newblock In Edward~N. Zalta and Uri Nodelman, editors, {\em The {Stanford} Encyclopedia of Philosophy}. Metaphysics Research Lab, Stanford University, {W}inter 2022 edition, 2022.

\bibitem{halbach2014axiomatic}
Volker Halbach.
\newblock {\em Axiomatic theories of truth}.
\newblock Cambridge University Press, 2014.

\bibitem{kunen2009foundations}
K.~Kunen.
\newblock {\em The Foundations of Mathematics}.
\newblock Mathematical logic and foundations. College Publications, 2009.

\bibitem{sep-prop-attitude-reports}
Michael Nelson.
\newblock {Propositional Attitude Reports}.
\newblock In Edward~N. Zalta and Uri Nodelman, editors, {\em The {Stanford} Encyclopedia of Philosophy}. Metaphysics Research Lab, Stanford University, {S}pring 2023 edition, 2023.

\bibitem{quine1956quantifiers}
Willard~V Quine.
\newblock Quantifiers and propositional attitudes.
\newblock {\em the Journal of Philosophy}, 53(5):177--187, 1956.

\bibitem{sep-logic-substructural}
Greg Restall.
\newblock {Substructural Logics}.
\newblock In Edward~N. Zalta, editor, {\em The {Stanford} Encyclopedia of Philosophy}. Metaphysics Research Lab, Stanford University, {S}pring 2018 edition, 2018.

\bibitem{sep-curry-paradox}
Lionel Shapiro and Jc~Beall.
\newblock {Curry’s Paradox}.
\newblock In Edward~N. Zalta, editor, {\em The {Stanford} Encyclopedia of Philosophy}. Metaphysics Research Lab, Stanford University, {W}inter 2021 edition, 2021.

\bibitem{sep-truth-values}
Yaroslav Shramko and Heinrich Wansing.
\newblock {Truth Values}.
\newblock In Edward~N. Zalta, editor, {\em The {Stanford} Encyclopedia of Philosophy}. Metaphysics Research Lab, Stanford University, {W}inter 2021 edition, 2021.

\end{thebibliography}
\end{document}